\documentclass{amsart}
\usepackage{lastpage}
\usepackage{yhmath,verbatim,mathabx}
\usepackage[all]{xy}

\DeclareFontFamily{U}{mathx}{\hyphenchar\font45}
\DeclareFontShape{U}{mathx}{m}{n}{
	<5> <6> <7> <8> <9> <10>
	<10.95> <12> <14.4> <17.28> <20.74> <24.88>
	mathx10
}{}
\DeclareSymbolFont{mathx}{U}{mathx}{m}{n}
\DeclareFontSubstitution{U}{mathx}{m}{n}
\DeclareMathAccent{\widecheck}{0}{mathx}{"71}
\DeclareMathAccent{\wideparen}{0}{mathx}{"75}

\newcommand{\bdis}{\begin{displaymath}}
\newcommand{\edis}{\end{displaymath}}
\newcommand{\be}{\begin{equation}}
\newcommand{\ee}{\end{equation}}
\newcommand{\mbb}{\mathbb}
\newcommand{\mcal}{\mathcal}

\newcommand{\vp}{\varphi}

\newcommand{\zf}{\zeta\left(\frac{1}{2}+it\right)}

\DeclareMathOperator{\dn}{dn}

\DeclareMathOperator{\sn}{sn}
\DeclareMathOperator{\cn}{cn}

\theoremstyle{definition}

\theoremstyle{remark}
\newtheorem{remark}[]{Remark}

\newtheorem*{mydef1}{{\bf Theorem}}

\newtheorem*{mydef5}{{\bf Lemma}}

\numberwithin{equation}{section}

\begin{document}

\title{Jacob's ladders and exact meta-functional equations on level curves as global quantitative characteristics of synergetic phenomenons excited by the function $|\zf|^2$}  

\author{Jan Moser}

\address{Department of Mathematical Analysis and Numerical Mathematics, Comenius University, Mlynska Dolina M105, 842 48 Bratislava, SLOVAKIA}

\email{jan.mozer@fmph.uniba.sk}


\begin{abstract}
In this paper we use operation of crossbreeding on the set of six transmutations of corresponding asymptotic complete hybrid formulas from our previous paper. We obtain in result the set of fifteen exact meta-functional equations. Every of them represents new formula in the theory of the Riemann's zeta-function. 
\begin{center} 
DEDICATED TO 160th ANNIVERSARY OF RIEMANN's FUNCTIONAL EQUATION.  
\end{center}
\end{abstract}
\maketitle 

\section{Introduction}

\subsection{} 

Let us remind that in the paper \cite{8} we have obtained the following: sets of values 
\be \label{1.1} 
\begin{split}
& \left\{\left|\zf\right|^2\right\},\ \{f_1(t)\}=\{\sin^2t\},\ \{f_2(t)\}=\{\cos^2t\},\\ 
& \{f_3(t)\}=\{\cos 2t\},\ t\in [\pi L,\pi L+U],\ U\in (0,\pi/4),\ L\in\mbb{N} 
\end{split}
\ee 
generates the following secondary asymptotic complete hybrid formula (see \cite{8}, (3.7), $k_1=k_2=k_3=1$) 
\be \label{1.2} 
\begin{split}
& \left|\zeta\left(\frac 12+i\alpha_1^{1,1}\right)\right|^2\sin^2\alpha_0^{1,1}-\\ 
& - \left\{1+\mcal{O}\left(\frac{\ln\ln L}{\ln L}\right)\right\}\left|\zeta\left(\frac 12+i\alpha_1^{2,1}\right)\right|^2\cos^2\alpha_0^{2,1}+\\ 
& + \left|\zeta\left(\frac 12+i\alpha_1^{3,1}\right)\right|^2\cos(2\alpha_0^{3,1})=0,\\ 
& \forall\- L\geq L_0<0, 
\end{split}
\ee 
($L_0$ beeing a sufficiently big constant), where 
\be \label{1.3} 
\begin{split}
& \alpha_r^{l,1}=\alpha_r(U,\pi L,f_l),\ r=0,1,\ l=1,2,3, \\ 
& \alpha_0^{l,1}\in (\pi L,\pi L+U),\ \alpha_1^{l,1}\in (\overset{1}{\wideparen{\pi L}},\overset{1}{\wideparen{\pi L+U}}), 
\end{split}
\ee  
and the mother formula (\ref{1.2}) generates the set of six transmutations. 

In this paper we use the operation of crossbreeding (see \cite{4} -- \cite{7}) on the set of these transmutations to obtain the set of exact meta-functional equations on level curves in the Gauss' plane. For example, there are sets 
\bdis 
\overset{(1)}{\Omega}_l,\ \overset{(5)}{\Omega}_l,\ \overset{(7)}{\Omega}_l(k_l)\subset\mbb{C},\ l=1,2,3 
\edis 
such that they fulfill following initial conditions 
\bdis 
\begin{split}
& \frac 12+i\alpha_1^{l,1}(U,\pi L,f_l)\in \overset{(1)}{\Omega}_l(U,L),\ l=1,2,3, \\ 
& \forall\- U\in(0,\pi/4),\ \forall\- L\geq L_0>0 
\end{split}
\edis 
such that for every of elements 
\bdis 
s_l^1\in \overset{(1)}{\Omega}_l,\ s_l^5\in \overset{(5)}{\Omega}_l,\ s_l^7\in \overset{(7)}{\Omega}_l 
\edis  
we have the following exact meta-functional equation 
\be \label{1.4} 
\begin{split}
& \frac{|\zeta(s_1^1)|^2}{|\Gamma(s_1^5)|^2}|\cn(s_2^7,k_2)|^2+\frac{|\zeta(s_3^1)|^2}{|\Gamma(s_3^5)|}|\cn(s_2^7,k_2)|^2=\\ 
& = \frac{|\zeta(s_1^1)|^2}{|\Gamma(s_2^5)|^2}|\sn(s_1^7,k_1)|^2+ \frac{|\zeta(s_3^1)|^2}{|\Gamma(s_2^5)|^2}|\dn(s_3^7,k_3)|
\end{split}
\ee 
for the corresponding subsets of the following sets 
\be \label{1.5} 
\begin{split}
& \{|\zeta(s)|\},\ \{|\Gamma(s)|\},\ \{\sn(s,k_1)\},\ \{|\cn(s,k_2)|\},\ \{|\dn(s,k_3)|\}, \\ 
& [(k_1)^2,(k_2)^2,(k_3)^2]\in (0,1)^3. 
\end{split}
\ee 

\begin{remark}
The exact meta-functional equation (\ref{1.4}) (together with the set of remaining 14 formulas) represents new type of formulas in the theory of Riemann's zeta-function as the last level in the sequence: 
\begin{itemize}
	\item[(a)] three exact $\zeta$-factorization formulas (see \cite{7}, Lemmas 1 --3) , 
	\item[(b)] asymptotic complete hybrid formula, 
	\item[(c)] six transmutation of the above, 
	\item[(d)] fifteen exact meta-functional equations. 
\end{itemize}
\end{remark} 

\begin{remark}
Interpretation of meta-functional equation (\ref{1.1}) within our $\zeta$-alchemy is that it gives the global quantitative characterization of synergetic phenomenons that lie in the cooperative interactions between corresponding subsets of the set (\ref{1.5}). 
\end{remark} 

\begin{remark}
This paper is also based on new notions and methods in the theory of Riemann's zeta-function we have introduced in our series of 49 papers concerning Jacob's ladders. These can be found in arXiv[math.CA] starting with the paper \cite{1}. 
\end{remark}

\subsection{} 

Furthermore, we give some remarks concerning connections between basic functions 
\bdis 
\left|\zf\right|^2,\ \vp_1(t),\ \vp_1^1(t),\ \alpha_r^{l,1},\ r=0,1,\ l=1,2,3. 
\edis 
Since the function $2\vp_1(T)$ is fixed solution to the integral equation (see \cite{1}, comp. \cite{8}, subsection 2.1) 
\bdis 
\int_0^{\mu[x(T)]}\left|\zf\right|^2e^{-\frac{2}{x(T)}t}{\rm d}t=\int_0^T\left|\zf\right|^2{\rm d}t, 
\edis  
then 
\be \label{1.6} 
\left|\zf\right|^2\longrightarrow \vp_1(t), 
\ee 
i.e. Jacob's ladder is generated by the function $\left|\zf\right|^2$. 

Next, we have 
\be \label{1.7} 
\alpha_r^{l,1}=\vp_1^{1-r}(d_l),\ d_l=d(U,\pi L,f_l),\ r=0,1,\ l=1,2,3 
\ee  
(see \cite{3}, (6.3), comp. \cite{8}, (2.6)), where 
\be \label{1.8} 
\vp_1(t)\longrightarrow \vp_1^1(t) 
\ee 
(of course, the first reverse iteration is generated by the Jacob's ladder). 

Now we have the following connection (see (\ref{1.6}) -- (\ref{1.8}))  
\be \label{1.9} 
\left|\zf\right|^2  \longrightarrow \vp_1(t) \longrightarrow \vp_1^1(t), 
\ee  
and, in more details (\ref{1.3}): 
\be \label{1.10} 
\begin{split}
& \alpha_r^{l,1}=\alpha_r(U,\pi L,f_l;[|\zeta_{0,5}|^2]),\ |\zeta_{0,5}|^2=\left|\zf\right|^2, \\ 
& r=0,1,\ l=1,2,3. 
\end{split}
\ee  

\begin{remark}
In our paper \cite{8} we have obtained seven sets of level curves, namely: 
\be \label{1.11} 
\begin{split}
& \overset{(1)}{\Omega}_l(U,L;[|\zeta_{0,5}|^2]),\ \overset{(2)}{\Omega}_l(U,L;[|\zeta_{0,5}|^2]),\ \overset{(3)}{\Omega}_l(U,L;[|\zeta_{0,5}|^2]),\\ 
& \overset{(4)}{\Omega}_l(U,L;n_l,[|\zeta_{0,5}|^2]),\ \overset{(5)}{\Omega}_l(U,L;[|\zeta_{0,5}|^2]),\ \overset{(6)}{\Omega}_l(U,L;p_l,[|\zeta_{0,5}|^2]),\\ 
& \overset{(7)}{\Omega}_l(U,L;k_l,[|\zeta_{0,5}|^2]),\ l=1,2,3, 
\end{split}
\ee  
i.e. the main founder of the class of sets (\ref{1.11}) is the function 
\bdis 
\left|\zf\right|^2 
\edis  
together with three elementary functions in (\ref{1.1}). 
\end{remark} 

\section{The structure of $\mcal{O}$-terms in (\ref{1.1})} 

\subsection{} 

Let us remind that the sets of values 
\be \label{2.1} 
\begin{split}
& \{f_1(t)\}=\{\sin^2t\},\ \{f_2(t)\}=\{\cos^2t\},\ 
 \{f_3(t)\}=\{\cos 2t\}, \\ 
&  t\in [\pi L,\pi L+U],\ U\in (0,\pi/4),\ L\in\mbb{N} 
\end{split}
\ee  
generate the following exact secondary complete hybrid formula (see \cite{7}, (3.6), $k_1=k_2=k_3=1$) 
\be \label{2.2} 
\begin{split}
& \tilde{Z}^2(\alpha_1^{1,1})\sin^2\alpha_0^{1,1}+\tilde{Z}^2(\alpha_1^{3,1})\cos(2\alpha_0^{3,1})=\tilde{Z}^2(\alpha_1^{2,1})\cos^2\alpha_0^{2,1}, \\ 
& L\geq L_0>0, 
\end{split}
\ee 
(where $L_0$ is sufficiently big), where 
\be \label{2.3} 
\tilde{Z}^2(t)=\frac{\left|\zf\right|^2}{\omega(t)},\ \omega(t)=\left\{1+\mcal{O}\left(\frac{\ln \ln t}{\ln t}\right)\right\}\ln t, 
\ee  
(see \cite{2}, (9.1), (9.2)). Next, we have 
\be \label{2.4} 
\begin{split}
& \alpha_r^{l,1}=\alpha_r(U,\pi L,f_l),\ r=0,1,\ l=1,2,3, \\ 
& \alpha_0^{l,1}\in (\pi L,\pi L+U),\ \alpha_1^{l,1}\in (\overset{1}{\wideparen{\pi L}},\overset{1}{\wideparen{\pi L+U}}), 
\end{split}
\ee 
and 
\bdis 
[\overset{1}{\wideparen{\pi L}},\overset{1}{\wideparen{\pi L+U}}]
\edis 
is the first reverse iteration (by means of the Jacob's ladder, see \cite{3}) of the basic segment 
\bdis 
[\pi L,\pi L+U]=[\overset{0}{\wideparen{\pi L}},\overset{0}{\wideparen{\pi L+U}}]. 
\edis 

\begin{remark}
The components of the main $\zeta$-disconncted set (for our case) 
\be \label{2.5} 
\Delta(\pi L,U,1)=[\pi L,\pi L+U]\bigcup [\overset{1}{\wideparen{\pi L}},\overset{1}{\wideparen{\pi L+U}}]
\ee 
are separated each from other by gigantic distance $\rho$, see \cite{3}, (5.12), comp. \cite{6}, (2.2) -- (2.9): 
\be \label{2.6} 
\rho\{[\pi L,\pi L+U],[\overset{1}{\wideparen{\pi L}},\overset{1}{\wideparen{\pi L+U}}]\}\sim 
\pi (1-c)\frac{L}{\ln L},\ L\to\infty, 
\ee 
where $c$ stands for the Euler's constant. 
\end{remark} 

\begin{remark}
Disconnected set (\ref{2.5}) has the following properties (see \cite{3}, (2.5 ) -- (2.7), $T=\pi L, r=1$): 
\begin{itemize} 
\item[(a)] lengths of its components are given by 
\be \label{2.7} 
|[\overset{r}{\wideparen{\pi L}},\overset{r}{\wideparen{\pi L+U}}]|=o\left(\frac{L}{\ln L}\right),\ r=0,1,\ L\to\infty, 
\ee  
\item[(b)] lengths of adjacent segments are given by 
\be \label{2.8} 
|[\pi L+U,\overset{1}{\wideparen{\pi L}}]|\sim \pi (1-c)\frac{L}{\ln L},\ L\to\infty, 
\ee 
\item[(c)] of course, 
\be \label{2.9} 
\pi L<\overset{1}{\wideparen{\pi L}},\ U=o\left(\frac{L}{\ln L}\right). 
\ee 
\end{itemize} 
\end{remark} 

\subsection{} 

Next, we have (see (\ref{2.2}), (\ref{2.3})) 
\be \label{2.10} 
\omega(\alpha_1^{l,1})=\left\{1+\mcal{O}_l\left(\frac{\ln \ln \alpha_1^{l,1}}{\ln \alpha_1^{l,1}}\right)\right\}\ln \alpha_1^{l,1},\ l=1,2,3. 
\ee 
Since (see (\ref{2.4})) 
\bdis 
\alpha_1^{l,1}\in (\overset{1}{\wideparen{\pi L}},\overset{1}{\wideparen{\pi L+U}}), 
\edis  
and (see (\ref{2.7}) -- (\ref{2.9}))  
\bdis 
\begin{split}
& \alpha_1^{l,1}-\pi L=\alpha_1^{l,1}-\overset{1}{\wideparen{\pi L}}+\overset{1}{\wideparen{\pi L}}-(\pi L+U)+\pi L+U-\pi L= \\ 
& = \alpha_1^{l,1}-\overset{1}{\wideparen{\pi L}}+[\overset{1}{\wideparen{\pi L}}-(\pi L+U)]+U, \\ 
& \alpha_1^{l,1}-\pi L<|[\overset{1}{\wideparen{\pi L}},\overset{1}{\wideparen{\pi L+U}}]|+ |[\pi L+U,\overset{1}{\wideparen{\pi L}}]|+U=\\ 
& = o\left(\frac{L}{\ln L}\right)+\mcal{O}\left(\frac{L}{\ln L}\right)+o\left(\frac{L}{\ln L}\right)=\mcal{O}\left(\frac{L}{\ln L}\right), 
\end{split}
\edis  
then 
\be \label{2.11} 
\alpha_1^{l,1}=\alpha_1^{l,1}-\pi L+\pi L=\pi L+\mcal{O}\left(\frac{L}{\ln L}\right), 
\ee  
and consequently 
\be \label{2.12} 
\ln \alpha_1^{l,1}=\left\{1+\mcal{O}\left(\frac{1}{\ln^2 L}\right)\right\}\ln \pi L,\ l=1,2,3. 
\ee  
Since the function 
\bdis 
\frac{\ln\ln t}{\ln t},\ t>e^e 
\edis  
is decreasing, then 
\be \label{2.13} 
\frac{\ln\ln \alpha_1^{l,1}}{\ln\alpha_1^{l,1}}<\frac{\ln\ln\pi L}{\ln\pi L}=\mcal{O}_l\left(\frac{\ln \ln L}{\ln L}\right),\ L\to\infty. 
\ee   
Now, we have from (\ref{2.10}) by (\ref{2.12}) and (\ref{2.13}) 
\be \label{2.14} 
\begin{split}
& \omega(\alpha_1^{l,1})=\left\{1+\mcal{O}_l\left(\frac{\ln \ln L}{\ln L}\right)\right\}\left\{1+\mcal{O}\left(\frac{1}{\ln^2 L}\right)\right\}\ln \pi L= \\ 
& = \left\{1+\mcal{O}_l\left(\frac{\ln \ln L}{\ln L}\right)\right\}\ln \pi L= \\ 
& = \{1+o_l(1)\}\ln\pi L,\ l=1,2,3. 
\end{split}
\ee  

\subsection{} 

Further, we have (see (\ref{2.2}), (\ref{2.14})) the following formula 
\be \label{2.15} 
\begin{split}
& \frac{1}{1+o_1(1)}\left|\zeta\left(\frac 12+i\alpha_1^{1,1}\right)\right|^2\sin^2\alpha_0^{1,1}+ \\ & \frac{1}{1+o_3(1)}\left|\zeta\left(\frac 12+i\alpha_1^{3,1}\right)\right|^2\cos(2\alpha_0^{3,1})= \\ 
& \frac{1}{1+o_2(1)}\left|\zeta\left(\frac 12+i\alpha_1^{2,1}\right)\right|^2\cos^2\alpha_0^{2,1}, 
\end{split}
\ee 
and from this, by making use of the following notations 
\be \label{2.16} 
\begin{split} 
& A_1=\left|\zeta\left(\frac 12+i\alpha_1^{1,1}\right)\right|^2\sin^2\alpha_0^{1,1}>0, \\ 
& A_2=\left|\zeta\left(\frac 12+i\alpha_1^{2,1}\right)\right|^2\cos^2\alpha_0^{2,1}>0, \\ 
& A_3=\left|\zeta\left(\frac 12+i\alpha_1^{3,1}\right)\right|^2\cos(2\alpha_0^{3,1})>0, 
\end{split} 
\ee 
we obtain 
\be \label{2.17} 
\frac{A_1}{1+o_1(1)}+\frac{A_3}{1+o_3(1)}=\frac{A_2}{1+o_2(1)}. 
\ee 
Since 
\be \label{2.18} 
\begin{split}
& \frac{A_1}{1+o_1(1)}+\frac{A_3}{1+o_3(1)}=\{1+o_4(1)\}A_1+\{1+o_6(1)\}A_3=\\ 
& = A_1+A_3+o_4(1)A_1+o_6(1)A_3=\\ 
& = (A_1+A_3)\left[1+\frac{A_1}{A_1+A_3}o_4(1)+\frac{A_3}{A_1+A_3}o_6(1)\right], 
\end{split}
\ee 
then 
\be \label{2.19} 
\begin{split}
& (A_1+A_3)\{1+o_7(1)\}=\{1+o_5(1)\}A_2, \\ 
& A_1+A_3=\frac{1+o_5(1)}{1+o_7(1)}A_2, 
\end{split}
\ee 
and (see (\ref{2.14})) 
\be \label{2.20} 
A_1+A_3=\{1+o_8(1)\}A_2,\ o_8(1)=\mcal{O}\left(\frac{\ln\ln L}{\ln L}\right), 
\ee 
i.e. we have the following 

\begin{mydef5}
The system of sets (\ref{2.1}) generates the following asymptotic secondary complete hybrid formula 
\be \label{2.21} 
\begin{split}
& \left|\zeta\left(\frac 12+i\alpha_1^{1,1}\right)\right|^2\sin^2\alpha_0^{1,1}-\left\{1+\mcal{O}\left(\frac{\ln\ln L}{\ln L}\right)\right\}\left|\zeta\left(\frac 12+i\alpha_1^{2,1}\right)\right|^2\cos^2\alpha_0^{2,1}\\ 
& + \left|\zeta\left(\frac 12+i\alpha_1^{3,1}\right)\right|^2\cos(2\alpha_0^{3,1})=0, \\ 
& \forall\- L\geq L_0>0,\ U\in (0,\pi/4), 
\end{split}
\ee  
where 
\be \label{2.22} 
\mcal{O}\left(\frac{\ln\ln L}{\ln L}\right)=F(\alpha_1^{1,1},\alpha_1^{2,1},\alpha_1^{3,1},\alpha_0^{1,1},\alpha_0^{3,1}), 
\ee 
that is $F$ does not depend on $\alpha_0^{2,1}$. 
\end{mydef5}

\section{List of the first generation of exact meta-functional equations} 

Let us remind that we have obtained in our paper \cite{8} the set 
\be \label{3.1} 
\{(3.10), (4.5), (4.10), (5.5), (5.10), (5.15)\}
\ee 
of six transmutations of the mother formula (\ref{2.21}). 

\begin{remark}
It is clear (see Section 2) that the factor 
\be \label{3.2} 
1+\mcal{O}\left(\frac{\ln\ln L}{\ln L}\right)
\ee 
which is contained in every element of the set (\ref{3.1}) is identical one. 
\end{remark} 

Consequently, we may apply the operation of crossbreeding (see \cite{4} -- \cite{7}, here the elimination of the function (\ref{3.2})) on every two different elements of the set (\ref{3.2}). 

\begin{remark}
Let the symbol 
\bdis 
(3.10) \times (4.5) \Rightarrow 
\edis 
stand for the phrase \emph{we obtain by crosbreeding of the transmutations (3.10) and (4.5)}. 
\end{remark} 

We obtain, as a result of crossbreedings on the set (\ref{3.1}), the following. 

\begin{mydef1}
There are the sets 
\bdis 
\begin{split}
& \overset{(1)}{\Omega}_l(U,L),\ \overset{(2)}{\Omega}_l(U,L),\ \overset{(3)}{\Omega}_l(U,L),\ \overset{(4)}{\Omega}_l(U,L,n_l), \\ 
& \overset{(5)}{\Omega}_l(U,L),\ \overset{(6)}{\Omega}_l(U,L,p_l),\ \overset{(7)}{\Omega}_l(U,L,k_l)\subset \mbb{C}, 
\end{split}
\edis 
where 
\bdis 
\begin{split}
& U\in (0,\pi/4),\ L\geq L_0>0,\ (n_1,n_2,n_3)\in\mbb{N}^3, \\ 
& (p_1,p_2,p_3)\in\mbb{Z}^3,\ [(k_1)^2,(k_2)^2,(k_3)^2]\in (0,1)^3, 
\end{split}
\edis  
and the sets $\overset{(1)}{\Omega}_l$ fulfill the following initial conditions 
\be \label{3.3} 
\begin{split}
& \frac 12+i\alpha_1^{l,1}(U,\pi L,f_l)\in \overset{(1)}{\Omega}_l(U,L),\ l=1,2,3, \\ 
& \forall\- U\in (0,\pi/4),\ \forall\- L\geq L_0>0,  
\end{split}
\ee 
(comp. \cite{8}, (3.3)) such that for each of the elements 
\bdis 
s_l^n\in \overset{(n)}{\Omega}_l,\ l=1,2,3,\ n=1,2,\dots,7
\edis 
we have the following set of fifteen exact meta-functional equations: 
\be \label{3.4} 
\begin{split}
& (3.10) \times (4.5) \Rightarrow \\ 
& |\zeta(s_1^1)|^2|\zeta(s_1^2)|^2|\cos s_2^3|^2+|\zeta(s_3^1)|^2|\zeta(s_3^2)||\cos s_2^3|^2=\\ 
& =|\zeta(s_1^1)|^2|\zeta(s_2^2)|^2|\cos s_1^3|^2+|\zeta(s_3^1)|^2|\zeta(s_2^2)|^2|\cos s_3^3|, 
\end{split}
\ee 
\be \label{3.5} 
\begin{split}
	& (3.10) \times (4.10) \Rightarrow \\ 
	& |\zeta(s_1^1)|^2|\zeta(s_1^2)|^2|s_2^4|^{2n_2}+|\zeta(s_3^1)|^2|\zeta(s_3^2)||s_2^4|^{2n_2}=\\ 
	& =|\zeta(s_1^1)|^2|\zeta(s_2^2)|^2|s_1^4|^{2n_1}+|\zeta(s_3^1)|^2|\zeta(s_2^2)|^2|s_3^4|^{n_3}, 
\end{split}
\ee  
\be \label{3.6} 
\begin{split}
	& (3.10) \times (5.5) \Rightarrow \\ 
	& |\zeta(s_1^1)|^2|\zeta(s_1^2)|^2|\Gamma(s_1^5)|^2|\Gamma(s_3^5)|+|\zeta(s_3^1)|^2|\zeta(s_3^2)||\Gamma(s_1^5)|^2|\Gamma(s_3^5)|=\\ 
	& =|\zeta(s_1^1)|^2|\zeta(s_2^2)|^2|\Gamma(s_2^5)|^2|\Gamma(s_3^5)|+|\zeta(s_3^1)|^2|\zeta(s_2^2)|^2|\Gamma(s_1^5)|^2|\Gamma(s_2^5)|^2, 
\end{split}
\ee 
\be \label{3.7} 
\begin{split}
	& (3.10) \times (5.10) \Rightarrow \\ 
	& |\zeta(s_1^1)|^2|\zeta(s_1^2)|^2|J_{p_2}(s_2^6)|^2+|\zeta(s_3^1)|^2|\zeta(s_3^2)||J_{p_2}(s_2^6)|^2=\\ 
	& =|\zeta(s_1^1)|^2|\zeta(s_2^2)|^2|J_{p_1}(s_1^6)|^2+|\zeta(s_3^1)|^2|\zeta(s_2^2)|^2|J_{p_3}(s_3^6)|, 
\end{split}
\ee  
\be \label{3.8} 
\begin{split}
	& (3.10) \times (5.15) \Rightarrow \\ 
	& |\zeta(s_1^1)|^2|\zeta(s_1^2)|^2|\cn(s_2^7,k_2)|^2+|\zeta(s_3^1)|^2|\zeta(s_3^2)||\cn(s_2^7,k_2)|^2=\\ 
	& =|\zeta(s_1^1)|^2|\zeta(s_2^2)|^2|\sn(s_1^7,k_1)|^2+|\zeta(s_3^1)|^2|\zeta(s_2^2)|^2|\dn(s_3^7,k_3)|, 
\end{split}
\ee 
\be \label{3.9} 
\begin{split}
	& (4.5) \times (4.10) \Rightarrow \\ 
	& |\zeta(s_1^1)|^2|s_2^4|^{2n_2}|\cos s_1^3|^2+|\zeta(s_3^1)|^2|s_2^4|^{2n_2}|\cos s_3^3|=\\ 
	& =|\zeta(s_1^1)|^2|s_1^4|^{2n_1}|\cos s_2^3|^2+|\zeta(s_3^1)|^2|s_3^4|^{n_3}|\cos s_2^3|^2, 
\end{split}
\ee 
\be \label{3.10} 
\begin{split}
& (4.5) \times (5.5) \Rightarrow \\ 
& \frac{|\zeta(s_1^1)|^2}{|\Gamma(s_2^5)|^2}|\cos s_1^3|^2+\frac{|\zeta(s_3^1)|^2}{|\Gamma(s_2^5)|^2}|\cos s_3^3|=\\ 
& = \frac{|\zeta(s_1^1)|^2}{|\Gamma(s_1^5)|^2}|\cos s_2^3|^2+\frac{|\zeta(s_3^1)|^2}{|\Gamma(s_3^5)|}|\cos s_2^3|^2, 
\end{split}
\ee 
\be \label{3.11} 
\begin{split}
	& (4.5) \times (5.10) \Rightarrow \\ 
	& |\zeta(s_1^1)|^2|J_{p_2}(s_2^6)|^2|\cos s_1^3|^2+|\zeta(s_3^1)|^2|J_{p_2}(s_2^6)|^2|\cos s_3^3|=\\ 
	& =|\zeta(s_1^1)|^2|J_{p_1}(s_1^6)|^2|\cos s_2^3|^2+|\zeta(s_3^1)|^2|J_{p_3}(s_3^6)||\cos s_2^3|^2, 
\end{split}
\ee 
\be \label{3.12} 
\begin{split}
& (4.5) \times (5.15) \Rightarrow \\ 
& |\zeta(s_1^1)|^2|\cos s_1^3|^2|\cn(s_2^7,k_2)|^2+|\zeta(s_3^1)|^2|\cos s_3^3||\cn(s_2^7,k_2)|^2=\\ 
& |\zeta(s_1^1)|^2|\cos s_2^3|^2|\sn(s_1^7,k_1)|^2+|\zeta(s_3^1)|^2|\cos s_2^3|^2|\dn(s_3^7,k_3)|, 
\end{split}
\ee 
\be \label{3.13} 
\begin{split}
	& (4.10) \times (5.5) \Rightarrow \\ 
	& \frac{|\zeta(s_1^1)|^2}{|\Gamma(s_2^5)|^2}|s_1^4|^{2n_1}+\frac{|\zeta(s_3^1)|^2}{|\Gamma(s_2^5)|^2}|s_3^4|^{n_3}=\\ 
	& = \frac{|\zeta(s_1^1)|^2}{|\Gamma(s_1^5)|^2}|s_2^4|^{2n_2}+\frac{|\zeta(s_3^1)|^2}{|\Gamma(s_3^5)|}|s_2^4|^{2n_2}, 
\end{split}
\ee 
\be \label{3.14} 
\begin{split}
&  (4.10) \times (5.10) \Rightarrow \\  
& |\zeta(s_1^1)|^2|J_{p_2}(s_2^6)|^2|s_1^4|^{2n_1}+|\zeta(s_3^1)|^2|J_{p_2}(s_2^6)|^2|s_3^4|^{n_3}=\\
& |\zeta(s_1^1)|^2|J_{p_1}(s_1^6)|^2|s_2^4|^{2n_2}+|\zeta(s_3^1)|^2|J_{p_3}(s_3^6)||s_2^4|^{2n_2}, 
\end{split}
\ee 
\be \label{3.15} 
\begin{split}
& (4.10) \times (5.15) \Rightarrow \\ 
& |\zeta(s_1^1)|^2|s_1^4|^{2n_1}|\cn(s_2^7,k_2)|^2+|\zeta(s_3^1)|^2|s_3^4|^{n_3}|\cn(s_2^7,k_2)|^2=\\
& |\zeta(s_1^1)|^2|s_2^4|^{2n_2}|\sn(s_1^7,k_1)|^2+|\zeta(s_3^1)|^2|s_2^4|^{2n_2}|\dn(s_3^7,k_3)|, 
\end{split}
\ee 
\be \label{3.16} 
\begin{split}
& (5.5) \times (5.10) \Rightarrow \\ 
& \frac{|\zeta(s_1^1)|^2}{|\Gamma(s_1^5)|^2}|J_{p_2}(s_2^6)|^2+ \frac{|\zeta(s_3^1)|^2}{|\Gamma(s_3^5)|}|J_{p_2}(s_2^6)|^2=\\ 
& \frac{|\zeta(s_1^1)|^2}{|\Gamma(s_2^5)|^2}|J_{p_1}(s_1^6)|^2+
\frac{|\zeta(s_3^1)|^2}{|\Gamma(s_2^5)|^2}|J_{p_3}(s_3^6)|, 
\end{split}
\ee 
\be \label{3.17} 
\begin{split}
& (5.5) \times (5.15) \Rightarrow \\ 
& \frac{|\zeta(s_1^1)|^2}{|\Gamma(s_1^5)|^2}|\cn(s_2^7,k_2)|^2+
\frac{|\zeta(s_3^1)|^2}{|\Gamma(s_3^5)|}|\cn(s_2^7,k_2)|^2=\\ 
& \frac{|\zeta(s_1^1)|^2}{|\Gamma(s_2^5)|^2}|\sn(s_1^7,k_1)|^2+ 
\frac{|\zeta(s_3^1)|^2}{|\Gamma(s_2^5)|^2}|\dn(s_3^7,k_3)|, 
\end{split}
\ee 
\be \label{3.18} 
\begin{split}
& (5.10) \times (5.15) \Rightarrow \\ 
& |\zeta(s_1^1)|^2|J_{p_1}(s_1^6)|^2|\cn(s_2^7,k_2)|^2+
|\zeta(s_3^1)|^2|J_{p_3}(s_3^6)||\cn(s_2^7,k_2)|^2=\\ 
& |\zeta(s_1^1)|^2|J_{p_2}(s_2^6)|^2|\sn(s_1^7,k_1)|^2+ 
|\zeta(s_3^1)|^2|J_{p_2}(s_2^6)|^2|\dn(s_3^7,k_3)|. 
\end{split}
\ee 
\end{mydef1} 

\begin{remark}
The set of initial conditions (\ref{3.3}) is an analogue of the Cauchy's initial conditions for a differential equation. These conditions, apart from other, eliminate some trivial manipulations with the level curves for an usual equations. 
\end{remark}

I would like to thank Michal Demetrian for his moral support of my study of Jacob's ladders.

\end{document}